\documentclass[3p,twocolumn]{elsarticle}

\usepackage{amsmath,amssymb,amsfonts,amsthm}
\usepackage{algorithm,algorithmicx}
\usepackage{algpseudocode}
\usepackage{graphicx}
\usepackage[hyphens]{url}
\usepackage{lineno,hyperref}
\hypersetup{colorlinks=true,breaklinks=true}
\usepackage{bbm}
\usepackage{mathtools}
\usepackage{csquotes}
\usepackage{enumitem}
\usepackage{cleveref}
\usepackage{subfigure}
\usepackage{booktabs}
\biboptions{sort&compress}
\modulolinenumbers[5]
\journal{Reliability Engineering \& System Safety }

\begin{document}

\begin{frontmatter}

\title{Near-optimal planning using approximate dynamic programming \\ to enhance post-hazard community resilience management}

\author[cive]{Saeed Nozhati\corref{cor1} \fnref{fn1}}
\ead{saeed.nozhati@colostate.edu}

\author[ece]{Yugandhar Sarkale\fnref{fn1}}
\ead{yugandhar.sarkale@colostate.edu}

\author[cive]{Bruce Ellingwood}
\ead{bruce.ellingwood@colostate.edu}

\author[ece,math]{Edwin K.P. Chong}
\ead{edwin.chong@colostate.edu}

\author[cive]{Hussam Mahmoud}
\ead{hussam.mahmoud@colostate.edu}

\cortext[cor1]{Corresponding author}
\fntext[fn1]{Contributed equally, listed alphabetically.}
\address[cive]{Dept. of Civil and Environmental Engineering, Colorado State University, Fort Collins, CO 80523-1372}
\address[ece]{Dept. of Electrical and Computer Engineering, Colorado State University, Fort Collins, CO 80523-1373}
\address[math]{Dept. of Mathematics, Colorado State University, Fort Collins, CO 80523-1874}
\begin{abstract}
The lack of a comprehensive decision-making approach at the community level is an important problem that warrants immediate attention. Network-level decision-making algorithms need to solve large-scale optimization problems that pose computational challenges. The complexity of the optimization problems increases when various sources of uncertainty are considered. This research introduces a sequential discrete optimization approach, as a decision-making framework at the community level for recovery management.  The proposed mathematical approach leverages approximate dynamic programming along with heuristics for the determination of recovery actions. Our methodology overcomes the curse of dimensionality and manages multi-state, large-scale infrastructure systems following disasters.
We also provide computational results showing that our methodology not only incorporates recovery policies of responsible public and private entities within the community but also substantially enhances the performance of their underlying strategies with limited resources. The methodology can be implemented efficiently to identify near-optimal recovery decisions following a severe earthquake based on multiple objectives for an electrical power network of a testbed community coarsely modeled after Gilroy, California, United States. The proposed optimization method supports risk-informed community decision makers within chaotic post-hazard circumstances.
\end{abstract}

\begin{keyword}
Approximate Dynamic Programming, Combinatorial Optimization, Community Resilience,  Electrical Power Network, Rollout Algorithm
\end{keyword}

\end{frontmatter}


\section{Introduction}

In the modern era, the functionality of infrastructure systems is of significant importance in providing continuous services to communities and in supporting their public health and safety. Natural and anthropogenic hazards pose significant challenges to infrastructure systems and cause undesirable system malfunctions and consequences. Past experiences show that these malfunctions are not always inevitable despite design strategies like increasing system redundancy and reliability \cite{nach}. Therefore, a sequential rational decision-making framework should enable malfunctioned systems to be restored in a timely manner after the hazards. Further, post-event stressors and chaotic circumstances, time limitations, budget and resource constraints, and complexities in the community recovery process, which are twinned with catastrophe, highlight the necessity for a comprehensive risk-informed decision-making framework  for recovery management at the community level. A comprehensive decision-making framework must take into account indirect and delayed consequences of decisions (also called the post-effect property of decisions), which requires foresight or planning. Such a comprehensive decision-making system must also be able to handle large-scale scheduling problems that encompass large combinatorial decision spaces to make the most rational plans at the community level.

Developing efficient computational methodologies for sequential decision-making problems has been a subject of significant interest \cite{kochenderfer,bertsekasbook}. In the context of civil engineering, several studies have utilized the framework of dynamic programming for management of bridges and pavement maintenance \cite{meidani,frangopol,jiang1,jiang2,Shaf}. Typical methodological formulations employ principles of dynamic programming that utilize state-action pairs. In this study, we develop a powerful and relatively unexplored methodological framework of formulating large infrastructure problems as \textit{string-actions}, which will be described in Section~\ref{5.1.2.}. Our formulation does not require an explicit state-space model; therefore, it is shielded against the common problem of state explosion when such methodologies are employed. The sequential decision-making methodology presented here not only manages network-level infrastructure but also considers the interconnectedness and cascading effects in the entire recovery process that have not been addressed in the past studies.

Dynamic programming formulations frequently suffer from the \emph{curse of dimensionality}. This problem is further aggravated when we have to deal with large combinatorial decision spaces characteristic of community recovery. Therefore, using approximation techniques in conjunction with the dynamic programming formalism is essential. There are several approximation techniques available in the literature \cite{Bertsekasneuro,kuhn,medury,busoniu}. Here, we use a promising class of approximation techniques called \emph{rollout} algorithms. We show how rollout algorithms blend naturally with our string-action formulation. Together, they form a robust tool to overcome some of the limitations faced in the application of dynamic programming techniques to massive real-world problems. The proposed approach is able to handle the curse of dimensionality in its analysis and management of multi-state, large-scale infrastructure systems and data sources. The proposed methodology is also able to consider and improve the current recovery policies of responsible public and private entities within the community.

Among infrastructure systems, electrical power networks (EPNs) are particularly critical insofar as the functionality of most other networks, and critical facilities depend on EPN functionality and management. Hence, the method is illustrated in an application to recovery management of the modeled EPN in Gilroy, California following a severe earthquake. The illustrative example shows how the proposed approach can be implemented efficiently to identify near-optimal recovery decisions. The computed near-optimal decisions restored the EPN of Gilroy in a timely manner, for residential buildings as well as main food retailers, as an example of critical facilities that need electricity to support public health in the aftermath of hazards.

The remainder of this study is structured as follows. In Section~\ref{sec:2}, we introduce the background of system resilience and the system modeling used in this study. In Section~\ref{case}, we introduce the case study used in this paper. In Section~\ref{hazard}, we describe the earthquake modeling, fragility, and restoration assessments.
In Section~\ref{Optprob}, we provide a mathematical formulation of our optimization problem. In Section~\ref{Optsol}, we describe the solution method to solve the optimization problem. In Section~\ref{Res}, we demonstrate the performance of the rollout algorithm with the string-action formulation through multiple simulations. In Section~\ref{Conc}, we present a brief conclusion of this research.

\section{System Resilience}\label{sec:2}
\begin{figure}
  \centering
      \includegraphics[width=.5\textwidth]{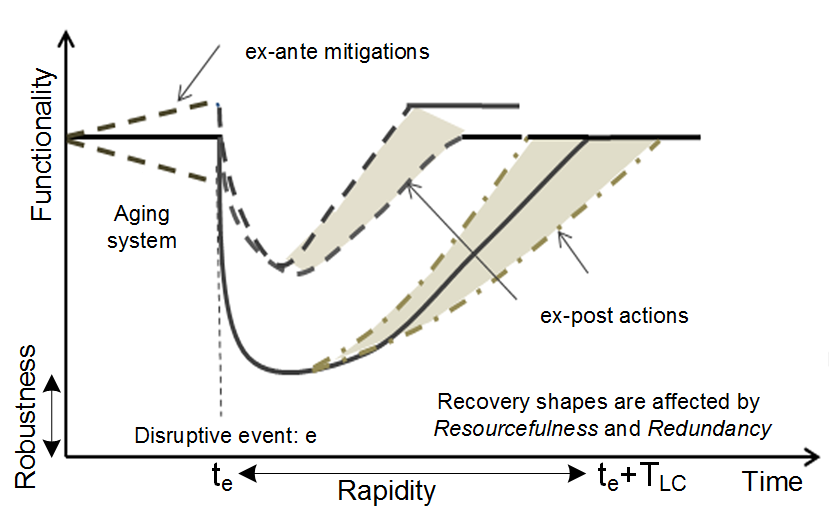}
  \caption{Schematic representation of
   resilience concept (adopted from \cite{mcallister,bruneau1})}
    \label{fig:1}
\end{figure}

The term \textit{resilience} is defined in a variety of ways. Generally speaking, resilience can be defined as \enquote{the ability to prepare for and adapt to changing conditions and withstand and recover rapidly from disruptions}\cite{directive2013critical}. Hence, resilience of a community (or a system) is usually delineated with the measure of community functionality, shown by the vertical axis of Fig.~\ref{fig:1} and four attributes of \textit{robustness}, \textit{rapidity}, \textit{redundancy}, and \textit{resourcefulness} \cite{bruneau1}. Fig.~\ref{fig:1} illustrates the concept of functionality, which can be defined as the ability of a system to support its planned mission, for example, by providing electricity to people and facilities. The understanding of interdependencies among the components of a system is essential to quantify system functionality and resilience. These interdependencies produce cascading failures where a large-scale cascade may be triggered by the malfunction of a single or few components \cite{buldyrev2010catastrophic}. Further, they contribute to the recovery rate and difficulty of the entire recovery process of a system. Different factors affect the recovery rate of a system, among which modification before disruptive events (ex-ante mitigations), different recovery policies (ex-post actions), and nature of the disruption are prominent \cite{barabadi2018post}. Fig.~\ref{fig:1} also highlights different sources of uncertainty that are associated with community functionality assessment and have remarkable impacts in different stages from prior to the event to the end of the recovery process. Therefore, any employed model to assess the recovery process should be able to consider the impacts of the influencing parameters.

In this study, the dependency of networks is modeled through an adjacency matrix $\bold{A}=[x_{ij}]$, where $x_{ij}\in[0,1]$ indicates the magnitude of dependency between components \textit{i} and \textit{j} \cite{watts1998collective}. In this general form, the adjacency matrix \textbf{A} can be a time-dependent stochastic matrix to capture the uncertainties in the dependencies and probable time-dependent variations.

According to the literature, the resilience index $\mathfrak{R}$ for each system is defined by the following equation \cite{bruneau1,cimellaro2014community}:
\begin{equation}\label{res}
\mathfrak{R} = \int _{t_e}^{t_e+T_{LC}} \displaystyle{\frac{Q(t)}{T_{LC}}} dt
\end{equation}

\noindent
where $Q(t)$ is the functionality of a system at time $\textit{t}$, $T_{LC}$ is the control time of the system, and $t_e$ is the time of occurrence of event \textit{e}, as shown in Fig.~\ref{fig:1}.

\section{Description of Case Study}\label{case}
\begin{figure}
  \centering
      \includegraphics[width=.5\textwidth]{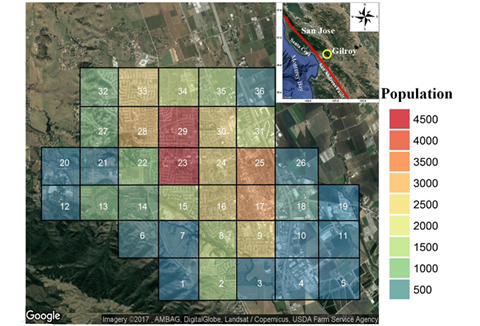}
  \caption{Map of Gilroy's population over the defined grids}
    \label{fig:2}
\end{figure}

In the case study of this paper, the community in Gilroy, California, USA is used as an example to illustrate the proposed approach. Gilroy is located approximately 50 kilometers (km) south of the city of San Jose with a population of 48,821 at the time of the 2010 census (see Fig.~\ref{fig:2}) \cite{Gilroy1}. The study area is divided into 36 gridded rectangles to define the community and encompasses 41.9 km\textsuperscript{2} area of Gilroy. In this study, we do not cover all the characteristics of Gilroy; however, the adopted model has a resolution that is sufficient to study the methodology at the community level under hazard events.

\begin{figure}
  \centering
      \includegraphics[width=.5\textwidth]{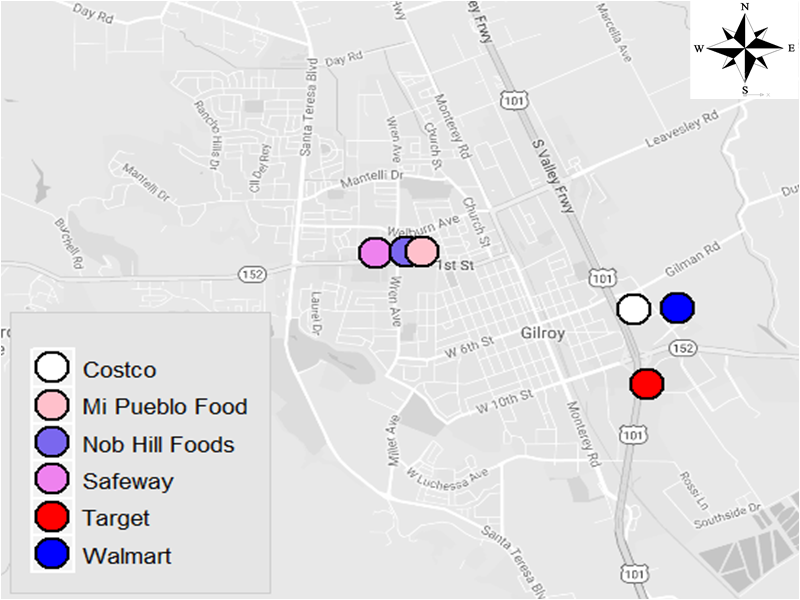}
  \caption{Gilroy's main food retailers}
    \label{fig:3}
\end{figure}

Gilroy contains 13 food retailers, but six main retailers, each of which has more than 100 employees, provide the main food requirements of Gilroy inhabitants \cite{Gilroy2}, as shown in  Fig.~\ref{fig:3} and summarized in Table \ref{T1}.

\begin{table*}[ht]
 \caption{The main food retailers of Gilroy }\label{T1}
\resizebox{\textwidth}{!}{
\begin{tabular}{lllllll}
   \hline
  Food Retailer & Walmart & Costco & Target & Mi Pueblo Food & Nob Hill Foods &Safeway\\
  \hline
Number of Employees & 395 & 220 & 130 & 106 & 100 & 130 \\
  \hline
\end{tabular}
}
\end{table*}

To assign the probabilities of shopping activity to each urban grid rectangle, the gravity model \cite{Popgen} is used. The gravity model identifies the shopping location probabilistically, given the location of residences. These probabilities are assigned to be proportional to food retailers\textsc{\char13} capacities and inversely corresponding to retailers\textsc{\char13} distances from centers of urban grid rectangles. Consequently, distant small locations are less likely to be selected than close large locations.

If the center of an urban grid is $c$, then food retailer $r$ is chosen according to the following distribution \cite{Popgen}:
\begin{equation}
P(r|c) \propto w_{r} e^{bT_{cr}}
\end{equation}
\noindent
where $w_{r}$ is the capacity of food retailer $r$, determined by Table 1, $b$ is a negative constant, and $T_{cr}$ is the travel time from urban grid rectangle $c$ to food retailer $r$. Google\textsc{\char13}s Distance Matrix API was called from within R by using the \textit{ggmap} package \cite{ggmap} to provide distances and travel times for the assumed transportation mode of driving.

\begin{figure}
  \centering
      \includegraphics[width=\linewidth]{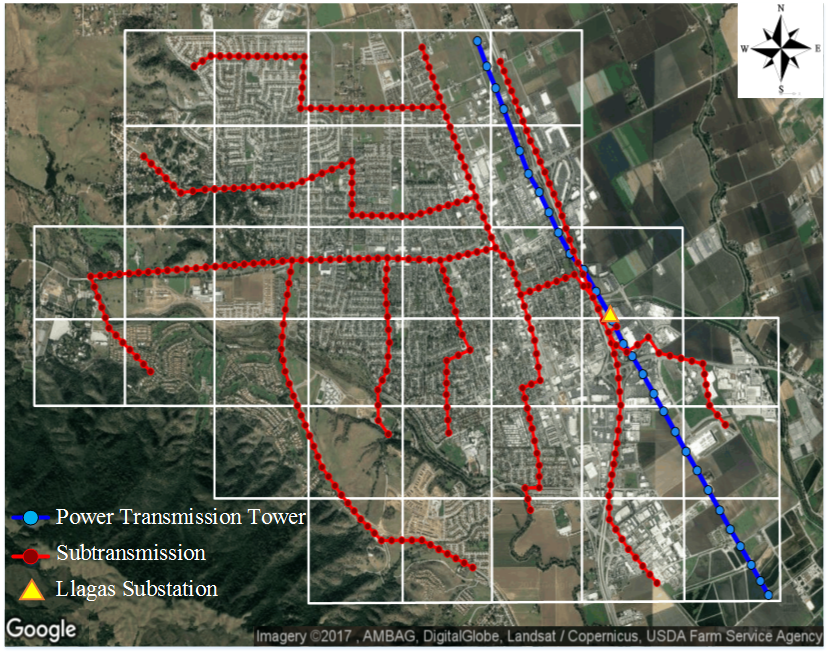}
  \caption{The modeled electrical power network of Gilroy}
    \label{fig:4}
\end{figure}

Fig.~\ref{fig:4} depicts the EPN components, located within the defined boundary. Llagas power substation, the main source of power in the defined boundary, is supplied by an 115~kV transmission line. Distribution line components are positioned at 100~m and modeled from the substation to the urban grids centers, food retailers, and water network facilities. In this study, the modeled EPN has 327 components.
\section{Hazard and Damage Assessment}\label{hazard}
  \subsection{Earthquake Simulation}
\begin{figure}
  \centering
      \includegraphics[width=\linewidth]{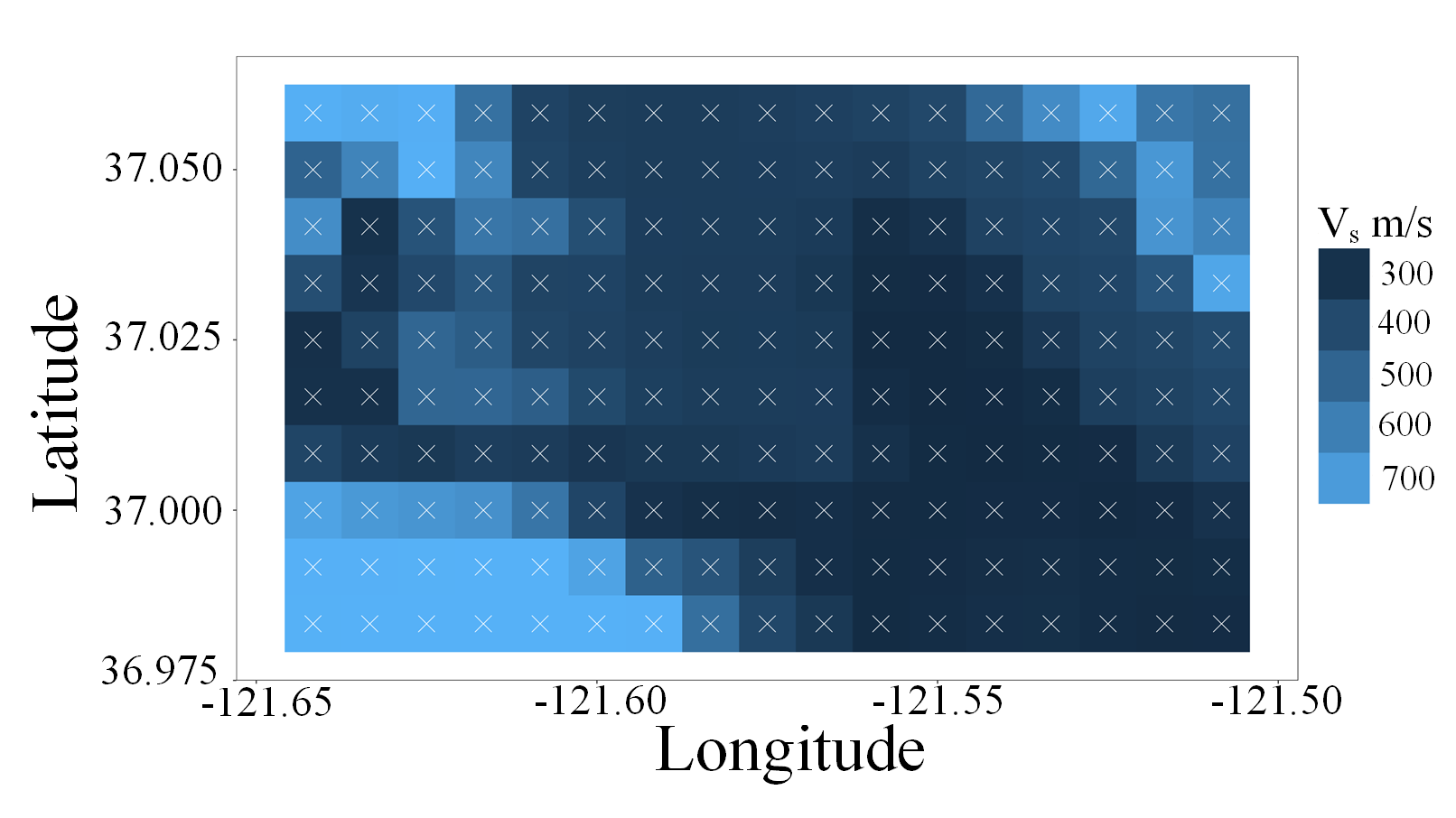}
  \caption{The map of shear velocity at Gilroy area}
    \label{fig:5}
\end{figure}

\begin{figure}
  \centering
      \includegraphics[width=\linewidth]{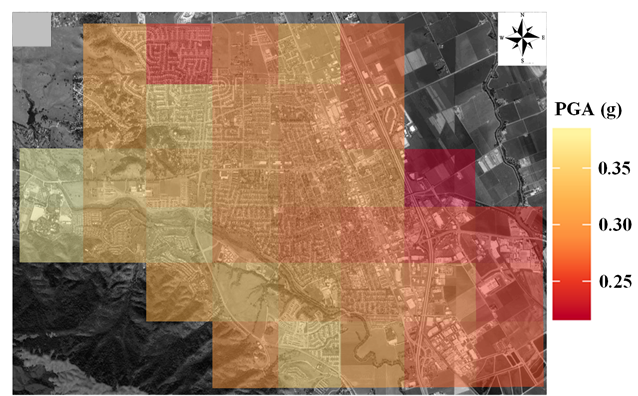}
  \caption{The simulation of median of peak ground acceleration field}
    \label{fig:6}
\end{figure}

  The seismicity of the Gilroy region of California is mainly controlled by the San Andreas Fault (SAF), which caused numerous destructive earthquakes like the Loma Prieta earthquake \cite{Loma1}. The spatial estimation of ground-motion amplitudes from earthquakes is an essential element of risk assessment, typically characterized by ground-motion prediction equations (GMPEs). GMPEs require several parameters, such as earthquake magnitude $M_{w}$, fault properties ($F_p$), soil conditions (i.e., the average shear-wave velocity in the top 30 m of soil, $V_{s30}$), and epicentral distances ($R$) to compute the seismic intensity measure ($\mathcal{IM}$) at any point. Modern GMPEs typically take the form
\begin{equation}
\begin{array}{r@{}l}
ln(\mathcal{IM}) &{}= f(M_w,R,V_{s30}, F_p) + \varepsilon_{1} \sigma + \varepsilon_{2} \tau\\[10pt]
ln(\mathcal{IM}) &{}= ln(\mathcal{\overline{IM}})+ \varepsilon_{1} \sigma + \varepsilon_{2} \tau
\end{array}
\end{equation}
where $\sigma$ and $\tau$ reflect the intra-event (within event) and inter-event (event-to-event) uncertainty respectively \cite{jayaram}. In this study, the GMPE proposed by Ambrahamson et al. \cite{abrahamson2013update} is used, and a ground motion similar to the Loma Prieta earthquake, one of the most devastating hazards that Gilroy has experienced \cite{Loma1}, with epicenter approximately 12 km of Gilroy downtown on the SAF projection is simulated. Figs.~\ref{fig:5} and \ref{fig:6} show the map of $V_{s30}$ and ground motion field for Peak Ground Acceleration (PGA), respectively.

  \subsection{Fragility Function and Restoration}
  In the event of an earthquake, the relations between ground-motion intensities and earthquake damage are pivotal elements in the loss estimation and the risk analysis of a community. Fragility curves describe the probability of experiencing or exceeding a particular level of damage as a function of hazard intensity.
   It is customary to model component fragilities with lognormal distributions \cite{kennedy1984seismic}. The conditional probability of being in or exceeding a particular damage state ($ds$), conditioned on a particular level of intensity measure $\mathcal{IM}=im$, is defined by
\begin{equation}
P(DS\geq ds | \mathcal{IM}=im)=\Phi\left(\frac{ln(im)-\lambda}{\xi}\right)
\end{equation}
where $\Phi$ is the standard normal distribution; $\lambda$ and $\xi$ are the mean and standard deviation of $ln(im)$.
The fragility curves can be obtained based on a) post-earthquake damage evaluation data (empirical curves) \cite{emp1} b) structural modeling (analytical curves) \cite{anl1} c) expert opinions (heuristics curves) \cite{hu1}. In the present study, the seismic fragility curves included in \cite{hazus, xie} are used for illustration.

To restore a network, a number of available resource units, $N$, as a generic single number including equipment, replacement components, and repair crews are considered for assignment to damaged components, and each damaged component is assumed to require only one unit of resource \cite{ouyang}. The restoration times based on the level of damage, used in this study, are presented in Table~\ref{T2}, based on \cite{hazus,saeedLA}.
\begin{table*}[ht]
  \caption{Restoration times based on the level of damage}\label{T2}
\resizebox{\linewidth}{!}{
\begin{tabular}{llllll}
  \hline
  &Damage States\\
   \hline
  Component & Undamaged & Minor & Moderate & Extensive & Complete \\
  \hline
  Electric sub-station &0& 1 & 3 & 7 & 30 \\
  Transmission line component&0& 0.5 & 1 & 1 & 2 \\
  Distribution line component &0& 0.5 & 1 & 1 & 1 \\
  \hline
\end{tabular}
}
\end{table*}\\

\section{Optimization Problem Description}\label{Optprob}
\subsection{Introduction}
 After an earthquake event occurs, each EPN component ends up in one of the damage states as shown in Table \ref{T2}. Let the total number of damaged components be $M$. Note that $M \leq 327 $. Both $M$ and $N$ are non-negative integers.  Also, in this study, $N\ll M$. This assumption is justified by the availability of limited resources with the planner where large number of components are damaged in the aftermath of a severe hazard.

 A decision maker or planner has the task of assigning units of resources to these damaged components. A decision maker has a heuristic or expert policy on the basis of which he can make his decisions to optimize multiple objectives. The precise nature of the objective of the planner can vary, which will be described in detail in Section~\ref{5.1.2.}. Particularly, at the first decision epoch, the decision maker or a resource planner deploys $N$ unit of resources at $N$ out of $M$ damaged components. Each unit of resource is assigned to a distinct damaged component. At every subsequent decision epoch, the planner must have an option of reassigning some or all of the resources to new locations based on his heuristics and objectives. He must have the flexibility of such a reassignment even if the repair work at the currently assigned locations is not complete. At every decision epoch, it is possible to forestall the reassignment of the units of resource that have not completed the repair work; however, we choose to solve the more general problem of preemptive assignment, where non-preemption at few or all the locations is a special case of our problem. The preemptive assignment problem is a richer decision problem than the non-preemptive case in the sense that the process of optimizing the decision actions is a more complex task because the size of the decision space is bigger.

In this study, we assume that the outcome of the decisions is fully predictable \cite{emi,iEMSs}. We are preparing a separate study to address the relaxation of this assumption, i.e., when outcomes of decisions exhibit uncertainty. The modified methods to deal with the stochastic problem will form a part of a separate paper \cite{ieeeconf}.

We improve upon the solutions offered by heuristics of the planner by formulating our optimization problem as a dynamic program and solving it using the rollout approach.

\subsection{Optimization Problem Formulation}\label{5.1.2.}
Suppose that the decision maker starts making decisions and assigning repair locations to different units of resource. The number of such non-trivial decisions to be made is less than or equal to $M-N$. When $M$ becomes less than or equal to $N$ (because of sequential application of repair actions to the damaged components), the assignment of units of resource becomes a trivial problem in our setting because each unit can simply be assigned one to one, in any order, to the damaged components. Consequently, a strict optimal assignment can be achieved in the trivial case. The size of this trivial assignment problem reduces by one for every new decision epoch until all the remaining damaged components are repaired. The additional units of resources retire because deploying more than one unit of resource to the same location does not decrease the repair time associated with that damaged component. Henceforth, we focus on the non-trivial assignment problem.

Let the variable $t$ denote the decision epoch, and let $D_t$ be the set of all damaged components before a repair action $x_t$ is performed. Let $t_{end}$ denote the decision epoch at which repair action $x_{t_{end}}$ is selected so that $|D_{t_{end}+1}|\leq N$. Note that $ t \in \mathcal{A} \coloneqq (1,2,\ldots,t_{end})$. Let $X=(x_1,x_2\ldots ,x_{t_{end}})$ represent the string of actions owing to the non-trivial assignment. We say that a repair action is completed when at least 1 out of the N damaged components is repaired.  Let $\mathcal{P}(D_t)$ be the powerset of $D_t$.
Let,
\begin{equation}
  \mathcal{P}_N(D_t)=\{C \in \mathcal{P}(D_t): |C|=N\}
\end{equation}
so that $x_t \in \mathcal{P}_N(D_t)$. Let $R_t$ be the set of all repaired components after the repair action $x_t$ is completed. Note that  $D_{t+1}=D_t\backslash R_t, ~ \forall t \in \mathcal{A}$, where $1 \leq |D_{t_{end}+1}| \leq N$, and the decision-making problem moves into the trivial assignment problem previously discussed.

We wish to calculate a string $X$ of repair actions that optimizes our objective functions $\mathit{F}(X)$. We deal with two objective functions in this study denoted by mapping $F_1$ and $F_2$.
\begin{itemize}
\item \textbf{Objective 1:} Let the variable $p$ represent the population of Gilroy and $\gamma$ represent a constant threshold. Let $X_1=(x_1,\ldots ,x_i)$ be the string of repair actions that results in restoration of electricity to $\gamma \times p$ number of people. Here, $x_i \in \mathcal{P}_N(D_i), $ where $D_i$ is the number of damaged component at the $i^{th}$ decision epoch. Let $n$ represent the time required to restore electricity to $\gamma \times p$ number of people as a result of repair actions $X_1$. 
    Formally,
    \begin{equation}\label{obj1}
      \mathit{F_1}(X_1)= n.
    \end{equation}
    Objective 1 is to compute the optimal solution $X_1^{*}$ given by
    \begin{equation}\label{obj1op}
      X_1^{*}=\arg\min_{X_1}  F_1(X_1).
    \end{equation}
    We explain the precise meaning of restoration of electricity to people in more detail in Section~\ref{Res_dis}. To sum up, in objective 1, our aim is to find a string of actions that minimizes the number of days needed to restore electricity to a certain fraction ($\gamma$) of the total population of Gilroy.
\item \textbf{Objective 2:} We define the mapping $F_2$ in terms of number of people who have electricity per unit of time; our objective is to maximize this mapping over a string of repair actions. Let the variable $k_t$ denote the total time elapsed between the completion of repair action $x_{t-1}$ and $x_t, ~ \forall  t \in \mathcal{A}\backslash \{1\}$, in which $k_1$ is the time elapsed between the start and completion of repair action $x_1$. Let $h_t$ be the total number of people who have electricity after the repair action $x_t$ is complete. 
    Then,
        \begin{equation}\label{obj2}
          F_2(X)=\frac{1}{k_{t_{end}}} \sum_{t=1}^{t_{end}}h_t \times k_t.
        \end{equation}
 We are interested in the optimal solution $X^*$ given by
            \begin{equation}\label{obj2op}
      X^{*}=\arg\max_{X}  F_2(X).
    \end{equation}
\end{itemize}

Note that our objective function in the second case $F_2(X)$ mimics the resilience index and can be interpreted in terms of \eqref{res}. Particularly, the integral in \eqref{res} is replaced by a sum because of discrete decision epochs, $Q(t)$ is replaced by the product $h_t \times k_t$, $k_{t_{end}}$ is analogous to $T_{LC}$ , and the integral limits are changed to represent the discrete decision epochs.

\section{Optimization Problem Solution}\label{Optsol}
Calculating $X^{*}$ or $X_1^{*}$ is a sequential optimization problem. The decision maker applies the repair action $x_t$ at the decision epoch $t$ to maximize or minimize a cumulative objective function. The string of actions, as represented in $X$ or $X_1$, are an outcome of this sequential decision-making process. This is particularly relevant in the context of dynamic programming where numerous solution techniques are available for the sequential optimization problem. Rollout is one such method that originated in dynamic programming. It is possible to use the dynamic programming formalism to describe the method of rollout, but here we accomplish this by starting from first principles \cite{Bertsekas2013}. We will draw comparisons between rollout with first principles and rollout in dynamic programming at suitable junctions. The description of the rollout algorithm is inherently tied with the notion of approximate dynamic programming.

\subsection{Approximate Dynamic Programming}\label{adp sec}
 Let's focus our attention on objective 1. The extension of this methodology to objective 2 is straightforward; we need to adapt notation used for objective 2 in the methodology presented below, and a maximization problem replaces a minimization problem. Recall that we are interested in the optimal solution  $X_1^*$ given by \eqref{obj1op}.  This can be calculated in the following manner: \\ First calculate $x^{*}_1$ as follows:
\begin{equation}\label{dp1}
x^{*}_1 \in \arg\min_{x_1} J_1(x_1)
\end{equation}
where the function $J_1$ is defined by
\begin{equation}\label{}
  J_1(x_1)=\min_{x_2,\ldots ,x_i}F_1(X_1).
\end{equation}
Next, calculate $x^{*}_2$ as:
\begin{equation}\label{}
  x^{*}_2 \in \arg\min_{x_2} J_2(x_1^*,x_2)
\end{equation}
where the function $J_2$ is defined by
\begin{equation}\label{}
  J_2(x_1,x_2)=\min_{x_3,\ldots ,x_i}F_1(X_1).
\end{equation}
Similarly, we calculate the $\alpha$-solution as follows:
\begin{equation}\label{som2}
  x^*_\alpha \in \arg\min_{x_\alpha}J_\alpha(x_1^*,\ldots ,x^*_{\alpha-1},x_\alpha)
\end{equation}
where the function $J_\alpha$ is defined by
\begin{equation}\label{som}
  J_\alpha(x_1,\ldots,x_\alpha)=\min_{x_{\alpha+1},\ldots ,x_i}F_1(X_1).
\end{equation}
The functions $J_\alpha$ are called the optimal cost-to-go functions and are defined by the following recursion:
\begin{equation}\label{}
  J_\alpha(x_1,\ldots,x_\alpha)=\min_{x_{\alpha+1}}J_{\alpha+1}(x_1,\ldots,x_{\alpha+1}),
\end{equation}
where the boundary condition is given by:
\begin{equation}\label{}
 J_i(X_1)=F_1(X_1).
\end{equation}
Note that $J$ is a standard notation used to represent cost-to-go functions in the dynamic programming literature.

The approach discussed above to calculate the optimal solutions is typical of the dynamic programming formulation. However, except for very special problems, such a formulation cannot be solved exactly because calculating and storing the \emph{optimal cost-to-go functions} $J_\alpha$ can be numerically intensive. Particularly, for our problem, let $|\mathcal{P}_N(D_t)|=\beta_t$; then the storage of $J_\alpha$ requires a table of size
\begin{equation}\label{tabeq}
  S_{\alpha}=\prod_{t=1}^{\alpha}\beta_{t},
\end{equation}
 where $\alpha \leq i$ for objective 1, and $\alpha \leq t_{end}$ for objective 2. In the dynamic programming literature, this is called as the \emph{curse of dimensionality}. If we consider objective 2 and wish to calculate $J_\alpha$ such that $\alpha=M-N$ (we assume for the sake of this example that only a single damaged component is repaired at each $t$), then for 50 damaged components and 10 unit of resources, $S_\alpha\approx 10^{280}$. In practice, $J_\alpha$ in \eqref{som2} is  replaced by an approximation denoted by $\tilde J_\alpha$. In the literature, $\tilde J_\alpha$ is called as a \emph{scoring function} or \emph{approximate cost-to-go} function \cite{Bertsekas1997}. One way to calculate $\tilde J_\alpha$ is with the aid of a heuristic; there are several ways to approximate $J_\alpha$ that do not utilize heuristic algorithms. All such approximation methods fall under the category of approximate dynamic programming.

 The method of rollout utilizes a heuristic in the approximation process. We provide a more detailed discussion on the heuristic in Section~\ref{Roll sec}. Suppose that a heuristic $\mathcal{H}$ is used to approximate the minimization in \eqref{som}, and let $H_\alpha(x_1,\ldots,x_\alpha)$ denote the corresponding approximate optimal value; then rollout yields the suboptimal solution by replacing $J_\alpha$ with $H_\alpha$ in \eqref{som2}:
\begin{equation}\label{}
  \tilde x_\alpha \in \arg\min_{x_\alpha}H_\alpha(\tilde x_1,\ldots ,\tilde x_{\alpha-1},x_\alpha).
\end{equation}

The heuristic used in the rollout algorithm is usually termed as the base heuristic. In many practical problems, rollout results in a significant improvement over the underlying base heuristic to solve the approximate dynamic programming problem \cite{Bertsekas1997}.

\subsection{Rollout Algorithm}\label{Roll sec}
It is possible to define the base heuristic $\mathcal{H}$ in several ways:
\begin{enumerate}[label=(\roman*)]
  \item The current recovery policy of regionally responsible public and private entities,
  \item The importance analyses that prioritize the importance of components based on the considered importance factors \cite{limnios},\label{o}
  \item The greedy algorithm that computes the greedy heuristic \cite{madant,liu},
  \item A random policy without any pre-assumption,
  \item A pre-defined empirical policy; e.g., base heuristic based on the maximum node and link \emph{betweenness} (shortest path), as for example, used in the studies of \cite{ouyang,masoomi}.\label{ll}
\end{enumerate}

The rollout method described in Section~\ref{adp sec}, using first principles and string-action formulation, for a discrete, deterministic, and sequential optimization problem has interpretations in terms of the policy iteration algorithm in dynamic programming. The policy iteration algorithm (see \cite{howard} for the details of the policy iteration algorithm including the definition of policy in the dynamic programming sense) computes an improved policy (policy improvement step), given a base policy (\emph{stationary}), by evaluating the performance of the base policy. The policy evaluation step is typically performed through simulations \cite{ieeeconf}. Rollout policy can be viewed as the improved policy calculated using the policy iteration algorithm after a single iteration of the policy improvement step. For a discrete and deterministic optimization problem, the base policy used in the policy iteration algorithm is equivalent to the base heuristic, and the rollout policy consists of the repeated application of this heuristic. This approach was used by the authors in \cite{Bertsekas1997} where they provide performance guarantees on the basic rollout approach and discuss variations to the rollout algorithm. Henceforth, for our purposes, base policy and base heuristic will be considered indistinguishable.

On a historical note, the term rollout was first coined by Tesauro in reference to creating computer programs that play backgammon \cite{Tesauro}. An approach similar to rollout was also shown much earlier in \cite{abramson}.

Ideally, we would like the rollout method to never perform worse than the underlying base heuristic (guarantee performance). This is possible under each of the following three cases \cite{Bertsekas1997}:
\begin{enumerate}
  \item The rollout method is \emph{terminating}  (called as optimized rollout). \label{case1}
  \item The rollout method utilizes a base heuristic that is \emph{sequentially consistent} (called as rollout). \label{case2}
  \item The rollout method is terminating and utilizes a base heuristic that is \emph{sequentially improving} (extended rollout and fortified rollout). \label{case3}
\end{enumerate}
  A sequentially consistent heuristic guarantees that the rollout method is terminating. It also guarantees that the base heuristic is sequentially improving. Therefore, \ref{case3} and \ref{case1} are the special cases of \ref{case2} with a less restrictive property imposed on the base heuristic (that of sequential improvement or termination). When the base heuristic is sequentially consistent, the fortified and extended rollout method are the same as the rollout method.

  A heuristic must posses the property of termination to be used as a base heuristic in the rollout method. Even if the base heuristic is terminating, the rollout method need not be terminating. Apart from the sequential consistency of the base heuristic, the rollout method is guaranteed to be terminating if it is applied on problems that exhibit special structure. Our problem exhibits such a structure. In particular, a finite number of damaged components in our problem are equivalent to the finite node set in \cite{Bertsekas1997}. Therefore, the rollout method in this study is terminating. In such a scenario, we could use the optimized rollout algorithm to guarantee performance without putting any restriction on the base heuristic to be used in the proposed formulation; however, a wiser base heuristic can potentially enhance further the computed rollout policy. Nevertheless, our problem does not require any special structure on the base heuristic for the rollout method to be sequentially improving, which is justified later in this section.

  In the terminology of dynamic programming, a base heuristic that admits sequential consistency is analogous to the Markov or stationary policy. Similarly, the terminating rollout method defines a rollout policy that is stationary.

Two different base heuristics are considered in this study. The first base heuristic is a random heuristic denoted by $\mathcal{H}$. The idea behind consideration of this heuristic is that in actuality there are cases where there is no thought-out strategy or the computation of such a scheme is computationally expensive. We will show though simulations that the rollout formulation can accept a random base policy at the community level from a decision maker and improve it significantly. The second base heuristic is called a \textit{smart} heuristic because it is based on the importance of components and expert judgment, denoted by $\mathcal{\hat H}$. The importance factors used in prioritizing the importance of the components can accommodate the contribution of each component in the network. This base heuristic is similar in spirit to the items \ref{o} and \ref{ll} listed above. More description on the assignment of units of resources based on $\mathcal{H}$ and $\mathcal{\hat H}$ is described in Section~\ref{sec_resl1}. We also argue there that $\mathcal{H}$ and $\mathcal{\hat H}$ are sequentially consistent. Therefore, in this study, and for our choice of heuristics, the extended, fortified, and rollout method are equivalent.

Let $\mathcal{H}$ be any heuristic algorithm; the state of this algorithm at the first decision epoch is $\tilde j_1$, where $\tilde j_1=(\tilde x_1)$. Similarly, the state of the algorithm at the $\alpha$\textsuperscript{th} decision epoch is the $\alpha$-solution given by $\tilde j_\alpha = (\tilde x_1,\ldots,\tilde x_\alpha)$, i.e., the algorithm generates the path of the states $(\tilde j_1,\tilde j_2,\ldots,\tilde j_\alpha)$. Note that $\tilde j_0$ is the dummy initial state of the algorithm $\mathcal H$. The algorithm $\mathcal{H}$ terminates when $\alpha=i$ for objective 1, and $\alpha=t_{end}$ for objective 2. Henceforth, in this section, we consider only objective 1 without any loss of generality. Let $H_\alpha(\tilde j_\alpha)$ denote the cost-to-go starting from the $\alpha$-solution, generated by applying $\mathcal{H}$ (i.e., $\mathcal{H}$ is used to evaluate the cost-to-go). The cost-to-go associated with the algorithm $\mathcal{H}$ is equal to the terminal reward, i.e., $\tilde H_\alpha(\tilde j_\alpha) = F_1(X_1)$. Therefore, we have: $\tilde H_1(\tilde j_1) = \tilde H_2(\tilde j_2) = \ldots = \tilde H_i(\tilde j_i)$. We use this heuristic cost-to-go in \eqref{som2} to find an approximate solution to our problem. This approximation algorithm is termed as ``Rollout on $\mathcal{H}$'' ($\mathcal{RH}$) owing to its structure that is similar to the approximate dynamic programming approach \emph{rollout}. The $\mathcal{RH}$ algorithm generates the path of the states $(j_1,j_2,\ldots,j_i)$ as follows:
\begin{equation}\label{rollmin}
j_{\alpha} = \arg\min_{\delta\in N(j_{\alpha-1})} \tilde J(\delta), ~~ \alpha=1,\ldots,i
\end{equation}
where, $j_{\alpha-1}=(x_1,\ldots,x_{\alpha-1})$, and
\begin{equation}\label{rollcal}
N(j_{\alpha-1}) = \{(x_1,\ldots,x_{\alpha-1},x) | x \in \mathcal{P}_N(D_\alpha)\}.
\end{equation}
The algorithm $\mathcal{RH}$ is sequentially improving with respect to $\mathcal{H}$ and outperforms $\mathcal{H}$ (see \cite{shankar2014} for the details of the proof).

The $\mathcal{RH}$ algorithm described above is termed as one-step lookahead approach because the repair action at any decision epoch $t$ (current step) is optimized by minimizing the cost-to-go given the repair action at $t$ (see \eqref{rollmin}). It is possible to generalize this approach to incorporate multi-step lookahead. Suppose that we optimize the repair actions at any decision epoch $t$ and $t+1$ (current and the next step combined) by minimizing the cost-to-go given the repair actions for the current and next steps. This can be viewed as a two-step lookahead approach. Note the similarity of this approach with the dynamic programming formulation from first principles in Section~\ref{adp sec}, except for the difficulty of estimating the cost-to-go values $J$ exactly. Also, note that a two-step lookahead approach is  computationally more intensive than the one-step approach. In principle, it is possible to extend it to step size $\lambda$, where $1\leq\lambda\leq i$. However, as $\lambda$ increases, the computational complexity of the algorithm increases exponentially. Particularly, when $\lambda$ is selected equal to $i$ at the first decision epoch, the $\mathcal{RH}$ algorithm finds the exact optimal solution by exhaustively searching through all possible combinations of repair action at each $t$, with computational complexity $O(S_i)$. Also, note that $\mathcal{RH}$ provides a tighter upper bound on the optimal objective value compared to the bound obtained from the original heuristic approach.

\section{Results}\label{Res}
\subsection{Discussion}\label{Res_dis}
We show simulation results for two different cases. In Case 1, we assume that people have electricity when their household units have electricity. Recall that the city is divided into different gridded rectangles according to population heat maps (Fig.~\ref{fig:2}), and different components of the EPN network serving these grids are depicted in Fig.~\ref{fig:4}. The entire population living in a particular gridded rectangle will not have electricity until all the EPN components serving that grid are either undamaged or repaired post-hazard (functional). Conversely, if the EPN components serving a particular gridded rectangle are functional, all household units in that gridded rectangle are assumed to have electricity.

In Case 2, along with household units, we incorporate food retailers into the analysis. We say that people have the \emph{benefit} of electric utility only when the EPN components serving their gridded rectangles are functional, and they go to a food retailer that is functional. A food retailer is functional (in the electric utility sense) when all the EPN components serving the retailer are functional. The mapping of number of people who access a particular food retailer is done at each urban grid rectangle and follows the gravity model explained in Section~\ref{case}.

In both the cases, the probability that a critical facility like a food retailer or an urban grid rectangle $G$ has electricity is
  \begin{equation}
    P(EG)\coloneqq P\left(\bigcap^{\hat n}_{l=1}EE_l\right)
  \end{equation}
   where $\hat n$ is the minimum number of EPN components required to supply electricity to $G$, $EG$ is the event that $G$ has electricity,  and $EE_l$ is the event that the $l^{th}$ EPN component has electricity. The sample space is a singleton set that has the outcome, \enquote{has electricity.}

For all the simulation results provided henceforth, the number of units of resource available with the planner is fixed at 10.

\subsubsection{Case 1: Repair Action Optimization of EPN for Household Units}\label{sec_resl1}
The search space $\mathcal{P}_N(D_t)$ undergoes a combinatorial explosion for modest values of $N$ and $D_t$, at each $t$, until few decision epochs before moving into the trivial assignment problem, where the value of $\beta_t$ is small. Because of the combinatorial nature of the assignment problem, we would like to reduce the search space for the rollout algorithm, at each $t$, without sacrificing on the performance. Because we consider EPN for only household units in this section, it is possible to illustrate techniques, to reduce the size of the search space for our rollout algorithm, that provide a good insight into formulating such methods for other similar problems. We present two representative methods to deal with the combinatorial explosion of the search space, namely, 1-step heuristic and N-step heuristic. Note that these heuristics are not the same as the base heuristic $\mathcal{H}$ or $\mathcal{\hat H}$.

Before we describe the 1-step and N-step heuristic, we digress to discuss $\mathcal{H}$ and $\mathcal{\hat H}$.
\begin{figure}
\includegraphics[width=\linewidth]{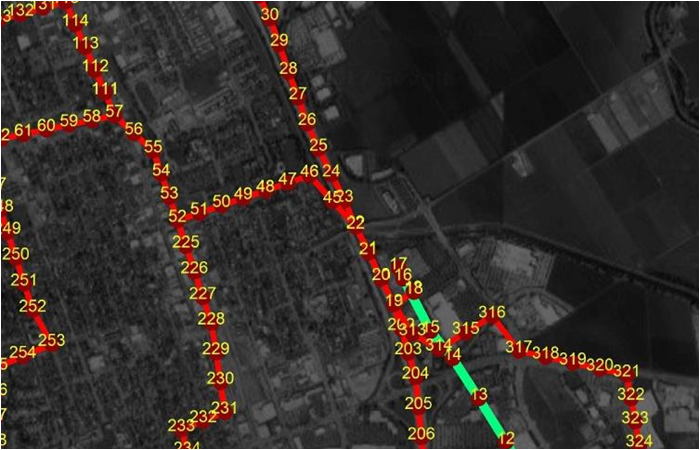}
\caption{Electrical Power Network of Gilroy}
\label{ZoomEPN}
\end{figure}
Both, $\mathcal{H}$ and $\mathcal{\hat H}$, have a preordained method of assigning units of resources to the damaged locations. This order of assignment remains fixed at each $t$. In $\mathcal{H}$, this order is decided randomly; while in $\mathcal{\hat H}$, it is decided based on importance factors. Let's illustrate this further with the help of an example. Suppose that we name each of the components of the EPN with serial numbers from 1 to 327 as shown partially in Fig.~\ref{ZoomEPN}; the assignment of these numbers to the EPN components is based on $\mathcal{\hat H}$ and remains fixed at each $t$, where a damaged component with a lower number is always assigned unit of resource before a damaged component with a higher number, based on the availability of units of resource. Therefore, the serial numbers depict the preordained priority assigned to components that is decided before decision-making starts. E.g., if the component number 21 and 22 are both damaged, the decision maker will assign one unit of resource to component 21 first and then schedule repair of component 22, contingent on availability of resources. Such a fixed pre-decided assignment of unit of resource by heuristic algorithm $\mathcal{H}$ and $\mathcal{\hat H}$ matches the definition of a consistent path generation in \cite{Bertsekas1997}. Therefore, $\mathcal{H}$ and $\mathcal{\hat H}$ are sequentially consistent. Note that the assignment of numbers 1 to 327 in Fig.~\ref{ZoomEPN} is assumed only for illustration purposes; the rollout method can incorporate a different preordained order defined by $\mathcal{H}$ and $\mathcal{\hat H}$.

\begin{figure}
\includegraphics[width=\linewidth]{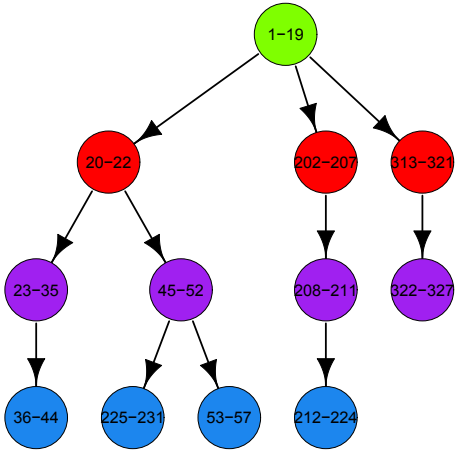}
\caption{Electrical Power Network Tree}
\label{depth}
\end{figure}

We now discuss the 1-step and N-step heuristic. In Fig.~\ref{ZoomEPN}, note that each successive EPN component (labeled 1-327), is dependent upon the prior EPN component for electricity. E.g., component 227 is dependent upon component 50. Similarly, component 55 is dependent upon component 50. The components in the branch 53-57 and 225-231 depend upon the component 52 for electricity. We exploit this serial nature of an EPN by representing the EPN network as a tree structure as shown in Fig.~\ref{depth}. Each number in the EPN tree represents an EPN component; each node represents a group of components determined by the label of the node, and the arcs of the tree capture the dependence of the nodes.

If the number of damaged components in the root node of the EPN tree is greater than $N$, then it would be unwise to assign a unit of resource at the first decision epoch to the fringe nodes of our EPN tree because we do not get any benefit until we repair damaged components in the root node. As soon as the number of damaged components in the root node of the EPN tree becomes less than $N$, only then we explore the assignment problem at other levels of the EPN tree.

\textbf{1-step Heuristic:}
We increase the pool of candidate damaged components, where the assignment of units of resources must be considered, to all the damaged components of the next level of the EPN tree if and only if the number of damaged components at the current level of the EPN tree is less than $N$. Even after considering the next level of the EPN tree, if the number of damaged components is less than $N$, we take one more step and account for all damaged components two levels below the current level. We repeat this until the pool of candidate damaged components is greater than or equal to $N$, or the levels of EPN tree are exhausted.

\textbf{N-step Heuristic:}
Note that it might be possible to ignore few nodes at each level of the EPN tree and assign units of resources to only some promising nodes. This is achieved in the N-step heuristic (here N in N-step is not same as $N$-number of workers). Specifically, if the number of the damaged components at the current level of the EPN tree is less than $N$, then the algorithm searches for a node at the next level that has the least number of damaged components, adds these damaged components to the pool of damaged components, and checks if the total number of damaged components at all explored levels is less than $N$. If the total number of candidate damaged components is still less than $N$, the previous process is repeated at unexplored nodes until the pool of damaged components is greater than or equal to $N$, or the levels of the EPN tree are exhausted.
Essentially, we do not consider the set ($D_t$) of all damaged components, at each $t$, but only a subset of $D_t$ denoted by $\tilde D^1_t$ (1-step heuristic) and $\tilde D^N_t$ (N-step heuristic).

We simulate multiple damage scenarios following an earthquake with the models discussed previously in Section~\ref{hazard}. On average, the total number of damaged components in any scenario exceeds 60\%.
\begin{figure}
\includegraphics[width=\linewidth]{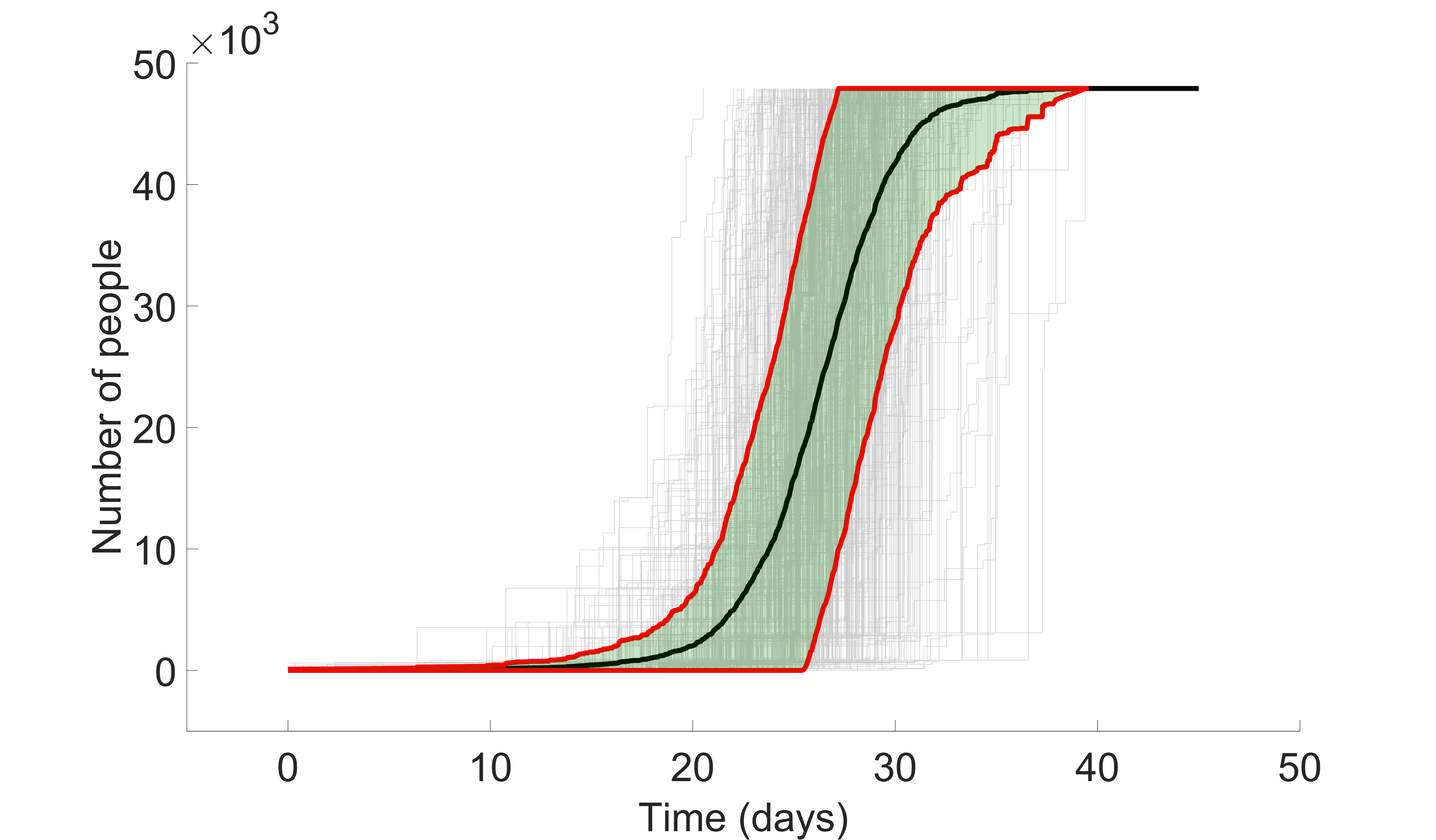}
\caption{Electrical power network recovery due to base heuristic $\mathcal{H}$ with one standard deviation band}
\label{Basenaive}
\end{figure}

We show the performance of $\mathcal H$ in Fig.~\ref{Basenaive}. The faint lines depict plots of EPN recovery for multiple scenarios when $\mathcal{H}$ is used for decision making. Here the objective pursued by the decision maker is objective 2. The black line shows the mean of all the recovery trajectories, and the red lines show the standard deviation. Henceforth, in various plots, instead of plotting the recovery trajectories for all the scenarios, we compare the mean of the different trajectories.
\begin{figure}
\includegraphics[width=\linewidth]{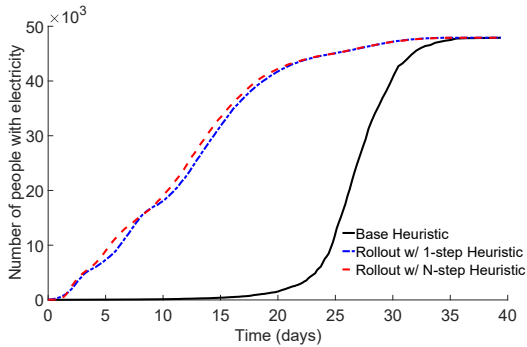}
\caption{Comparison of base heuristic ($\mathcal{H}$) vs. rollout algorithm with 1-step heuristic vs. rollout algorithm with N-step heuristic for multiple scenarios for objective 2.}
\label{compare}
\end{figure}
\begin{figure}[h]
\includegraphics[width=\linewidth]{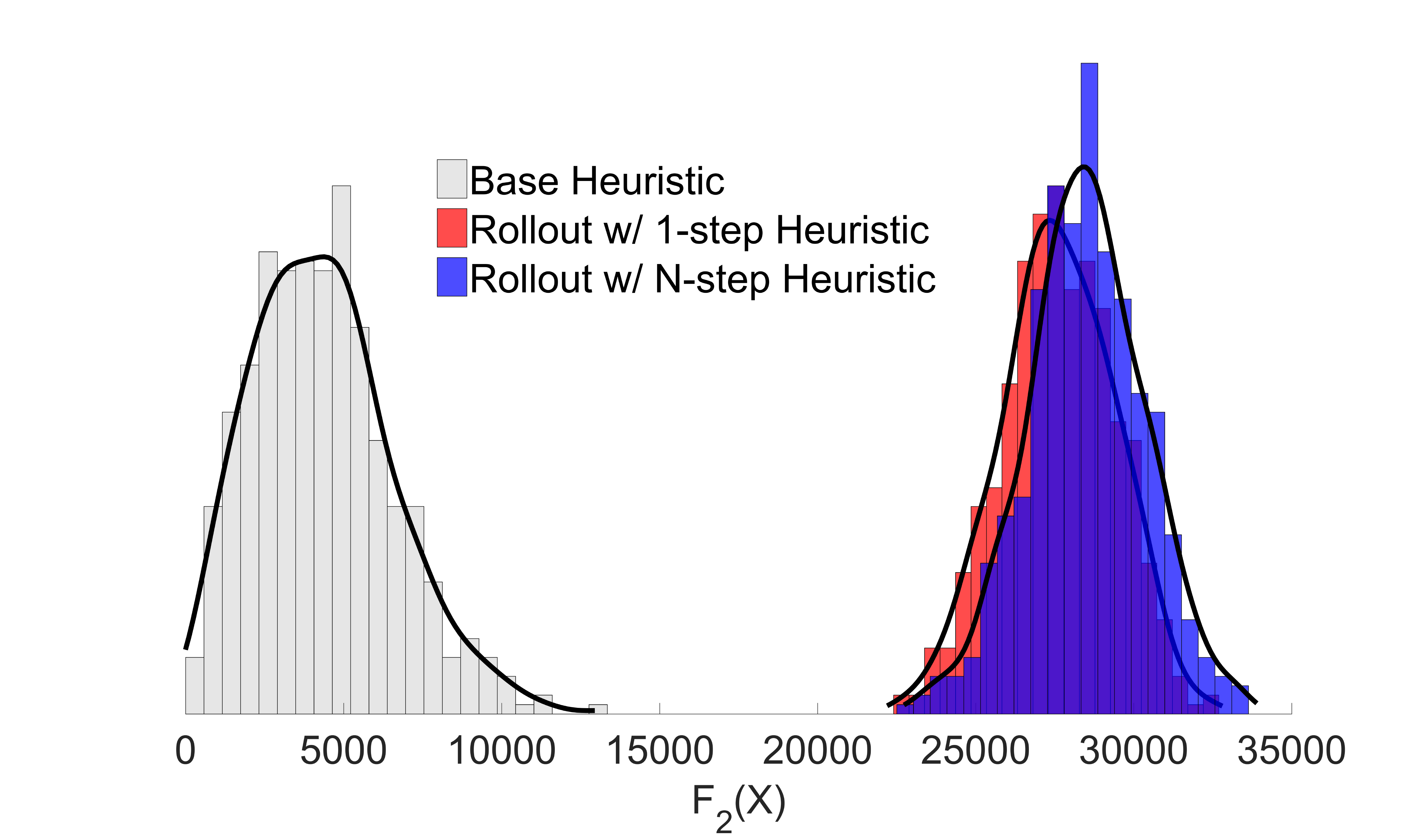}
\caption{Histogram of $F_2(X)$ with Base ($\mathcal{H}$), rollout with 1-step and  rollout with N-step heuristic for multiple scenarios}
\label{Qhis}
\end{figure}

Fig.~\ref{compare} shows the performance of our $\mathcal{RH}$ algorithm with respect to $\mathcal{H}$. Our simulation results demonstrate significant improvement over algorithm $\mathcal{H}$ when $\mathcal{RH}$ is used, for both the 1-step case and the N-step case. Another result is the performance shown by the 1-step heuristic with respect to the N-step heuristic. Even though the N-step heuristic skips some of the nodes at each level in EPN tree, to accommodate more promising nodes into the search process, the performance improvement shown is minimal. Even though all the damaged components are not used to define the search space of the rollout algorithm, only a small subset is chosen with the use of either 1-step and N-step heuristic (limited EPN tree search), the improvement shown by $\mathcal{RH}$ over $\mathcal{H}$ is significant. This is because pruning the search space of rollout algorithm using a subset of $D_t$ (restricting an exhaustive search), is only a small part of the entire rollout algorithm. Further explanation for such a behavior is suitably explained in Section~\ref{subcase2}.

Fig.~\ref{Qhis} shows the histogram of values of $F_2(X)$ for multiple scenarios, as a result of application of string-actions computed using $\mathcal{H}$ and $\mathcal{RH}$ (1-step and N-step heuristic). The rollout algorithms show substantial improvement over $\mathcal{H}$, for our problem.
\begin{figure}
\includegraphics[width=\linewidth]{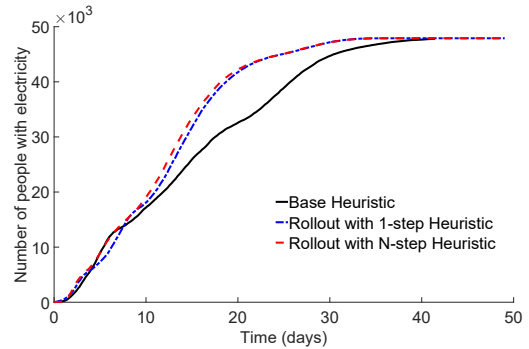}
\caption{Comparison of base heuristic ($\mathcal{\hat H}$) vs. rollout algorithm with 1-step heuristic vs. rollout algorithm with N-step heuristic for multiple scenarios for objective 2.}
\label{obj2smart}
\end{figure}

Fig.~\ref{obj2smart} shows simulation results for multiple scenarios when objective 2 is optimized, but when $\mathcal{\hat H}$ is considered instead of $\mathcal{H}$. This simulation study highlights interesting behavior exhibited by the rollout algorithm. In the initial phase of decision making, the algorithm $\mathcal{\hat H}$ competes with the rollout algorithm $\mathcal{R\hat H}$, slightly outperforming the rollout algorithm in many instances. However, after a period of 10 days, rollout (both 1-step and N-step heuristic) comprehensively outperforms $\mathcal{\hat H}$. Because rollout has the lookahead property, it offers conservative repair decisions initially (despite staying competitive with $\mathcal{\hat H}$), in anticipation of overcoming the loss suffered due to initial conservative decisions. Optimizing with foresight is an emergent behavior exhibited by our optimization methodology, which can offer significant advantages in critical decision-making scenarios.

\subsubsection{Case 2: Repair Action Optimization of EPN for Household Units \& Retailers}\label{subcase2}
\begin{figure}
\includegraphics[width=\linewidth]{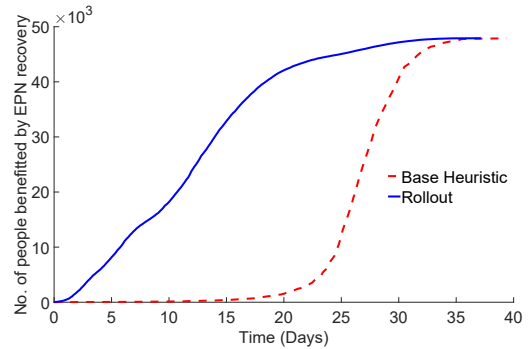}
\caption{Comparison of base heuristic ($\mathcal{H}$) vs rollout algorithm}
\label{detobj2}
\end{figure}

\begin{figure}
\includegraphics[width=\linewidth]{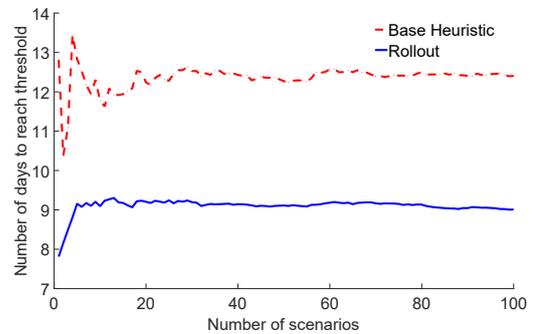}
\caption{Comparison of base heuristic ($\mathcal{\hat H}$)  vs. rollout algorithm for multiple scenarios for objective 1.}
\label{detobj1smart}
\end{figure}

\begin{figure}
\includegraphics[width=\linewidth]{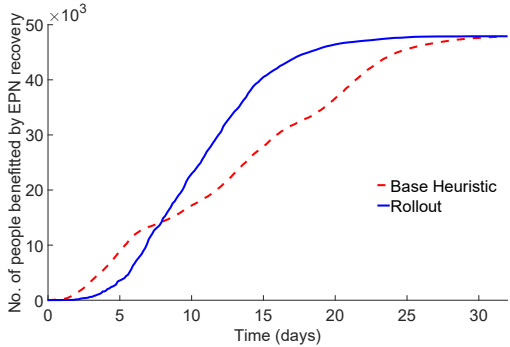}
\caption{Comparison of base heuristic ($\mathcal{\hat H}$) vs. rollout algorithm for multiple scenarios for objective 2.}
\label{detobj2smart}
\end{figure}

It is difficult to come up with techniques similar to 1-step and N-step heuristic to reduce the size of the search space $\mathcal{P}_N(D_t)$ for objective 1 and objective 2, when multiple networks are considered simultaneously in the analysis. This is because any such technique would have to simultaneously balance the pruning of candidate damaged components in Fig.~\ref{depth} (to form subsets like $\tilde D^1_t$ and $\tilde D^N_t$), serving both food retailers and household units. There is no natural/simple way of achieving this as in the case of the 1-step and N-step heuristics where only household units were considered. Whenever we consider complex objective functions and interaction between networks, it is difficult to prune the action space in a physically meaningful way just by applying heuristics.  For our case, any such heuristic will have to incorporate the gravity model explained in Section~\ref{case}. The heuristic must also consider the network topology and actual physical position of important EPN components within the network. As previously seen, our methodology works well even if we select only a small subset of $D_t$ ($\tilde D^1_t$ and $\tilde D^N_t$) to construct $\mathcal{P}_N(D_t)$ to avoid huge computational costs. This is because our methodology leverages the use of the one-step lookahead property with consistent and sequential improvement over algorithm $\mathcal{H}$ or $\mathcal{\hat H}$, to overcome any degradation in performance as a result of choosing $\tilde D_t \ll D_t$.
\begin{figure}
\includegraphics[width=\linewidth]{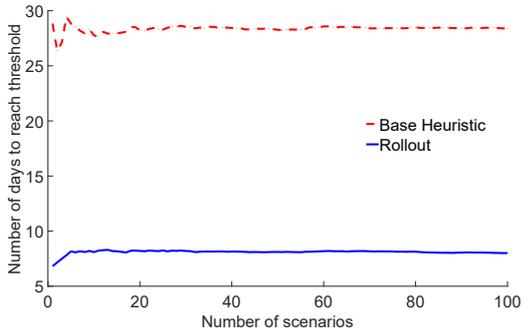}
\caption{Cumulative moving average plot for objective 1 where $\gamma=.8$ with base heuristic ($\mathcal{H}$) and rollout algorithm }
\label{detobj1}
\end{figure}
This is further justified in the simulation results shown in Figs.~\ref{detobj1} to \ref{detobj2smart}.

In Figs.~\ref{detobj1} and \ref{detobj2}, $\mathcal{H}$ is used as the base heuristic. This base heuristic is the same as the one used in the simulations shown in Case 1, which is not particularly well tuned for Case 2. Despite this, $\mathcal{RH}$ shows a stark improvement over $\mathcal{H}$. Fig.~\ref{detobj1} shows that the rollout algorithm, with the random selection of candidate damaged components ($\tilde D_t$), significantly outperforms $\mathcal{H}$ for objective 1. When we select the candidate damaged locations randomly, in addition to the randomly selected damaged components, we add to the set $\mathcal{P}_N(\tilde D_t)$ the damaged components selected by $\mathcal{\hat H}$ at each $t$.  For the rollout algorithm, the mean number of days, over multiple damage scenarios, to provide electricity to $\gamma$ times the total number of people is approximately 8 days, whereas for the base heuristic it is approximately 30. Similarly, Fig.~\ref{detobj2} shows that for objective 2 the benefit of EPN recovery per unit of time with rollout is significantly better than $\mathcal{H}$.

Figs.~\ref{detobj1smart} and \ref{detobj2smart} provide the synopsis of the results for both the objectives when $\mathcal{\hat H}$ is considered. In Fig.~\ref{detobj1smart}, the number of days to reach a threshold $\gamma=0.8$ as a result of $\mathcal{\hat H}$ algorithm is better than $\mathcal{H}$. However, note that the number of days to achieve objective 1 is still fewer using the rollout algorithm. The key inference from this observation is that rollout might not always significantly outperform the base heuristic but will never perform worse than the underlying heuristic.

For simulations in Figs.~\ref{detobj1smart} and \ref{detobj2smart}, the candidate damaged locations in defining the search space for the rollout algorithm are again chosen randomly and are a subset of $D_t$. As in the case of simulations in Fig.~\ref{detobj1}, we add damaged components selected by $\mathcal{\hat H}$ to the set $\mathcal{P}_N(\tilde D_t)$. Note that the number of days required to restore electricity to 80\% of people in Fig.~\ref{detobj1} is a day less than that required in Fig.~\ref{detobj1smart} despite performing rollout on a random base heuristic instead of the smart base heuristic used in the later. This can be attributed to 3 reasons: a) The damaged components in the set $(\tilde D_t)$ are chosen randomly for each simulation case. b) $\mathcal{\hat H}$ was designed for the simulations in Section.~\ref{sec_resl1} and is not particularly well tuned for simulations when both household units and retailers are considered simultaneously. c) $\mathcal{\hat H}$ is used in simulations in Fig.~\ref{detobj1smart} to approximate the cost-to-go function whereas $\mathcal{H}$ is used in simulations in Fig.~\ref{detobj1} for the approximation of the scoring function.

Fig.~\ref{detobj2smart} shows a behaviour similar to Fig.~\ref{obj2smart} (Case 1) where $\mathcal{\hat H}$  might outperform rollout in the short-term (as a result of myopic decision making on the part of the heuristic), but in the long run, rollout compressively improves upon the string-actions provided by algorithm $\mathcal{\hat H}$.

\subsection{Computational Efforts}
We provide a brief description of the computational efforts undertaken to optimize our decision-making problem. We have simulated multiple damage scenarios for each simulation result provided in this work. The solution methodology was implemented in MATLAB. The rollout algorithm gives multiple calls to the base heuristic function and searches for string-actions over the set $\mathcal{P}_N(\tilde D_t)$, $\mathcal{P}_N(\tilde D^1_t)$, or $\mathcal{P}_N(\tilde D^N_t)$. This is akin to giving millions of calls to the simulator. Therefore, our implementation of the code had to achieve a run time on the order of few micro seconds ($\mu s$) so that millions of calls to the simulator are possible. A single call to the simulator will be defined as evaluating some repair action $x_t \in \mathcal{P}_N(\tilde D_t)$, $\mathcal{P}_N(\tilde D^1_t)$, or $\mathcal{P}_N(\tilde D^N_t)$ and completing the rollout process until all the damaged components are repaired (complete rollout \cite{Bertsekasconst}). Despite the use of best software practices, to mitigate large action spaces by coding a fast simulator, it is imperative to match it with good hardware capability for simulations of significant size. It is possible to parallelize the simulation process to fully exploit modern day hardware capabilities.  We ran our simulations on the Summit super computer (see Section~\ref{ack}). Specifically, (100/376) Poweredge C6320 Haswell nodes were used, where each node has 2x e5-2680v3 (2400/9024) cores. In Case 1, we never encountered any $|\mathcal{P}_N(\tilde D^1_t)|$ or $|\mathcal{P}_N(\tilde D^N_t)|$ exceeding $10^6$ for any damage scenario, whereas for case 2 we limited the $|\mathcal{P}_N(\tilde D_t)|$ to at most $10^5$, for all damage scenarios. On this system, with a simulator call run time of 1~$\mu s$, we managed to run $10^{11}$ simulator calls per node for different simulation plots. Our computational efforts emphasize that for massive combinatorial decision-making problems, it is possible to achieve near-optimal performance with our methodology by harnessing powerful present day computational systems and implementing sound software practices.

\section{Concluding Remarks}\label{Conc}
In this study, we proposed an optimization formulation based on the method of rollout, for the identification of near-optimal community recovery actions, following a disaster. Our methodology utilizes approximate dynamic programming algorithms along with heuristics to overcome the curse of dimensionality in its analysis and management of large-scale infrastructure systems.
We have shown that the proposed approach can be implemented efficiently to identify near-optimal recovery decisions following a severe earthquake for the electrical power serving Gilroy. Different base heuristics, used at the community-level recovery, are used in the simulation studies. The intended methodology and the rollout policies significantly enhanced these base heuristics. Further, the efficient performance of the rollout formulation in optimizing different common objective functions for community resilience is demonstrated. We believe that the methodology is adaptable to other infrastructure systems and hazards.
\section{Acknowledgement}\label{ack}
“The research herein has been funded by the National Science Foundation under Grant CMMI-1638284.  This support is gratefully acknowledged.   Any opinions, findings, conclusions, or recommendations presented in this material are solely those of the authors and do not necessarily reflect the views of the National Science Foundation.”

“This work utilized the RMACC Summit supercomputer, which is supported by the National Science Foundation (awards ACI-1532235 and ACI-1532236), the University of Colorado Boulder and Colorado State University. The RMACC Summit supercomputer is a joint effort of the University of Colorado Boulder and Colorado State University.”

\section{References}

\bibliography{mybibfile}

\end{document}